\documentclass[a4paper,12pt]{article}

\usepackage{amsmath}
\usepackage{amssymb}
\usepackage{amsfonts}
\usepackage{graphicx}
\usepackage{mathrsfs}
\usepackage{pgfplots}

\newtheorem{theorem}{Theorem}

\title{\bf SPECTRA - a Maple library for solving linear matrix inequalities in exact arithmetic\footnote{This work was partly funded by the ERC Advanced Grant Taming.}}

\begin{document}

\author{
Didier Henrion\thanks{LAAS-CNRS, Universit\'e de Toulouse, CNRS, Toulouse, France;
Faculty of Electrical Engineering, Czech Technical University in Prague, Czech Republic.} \and 
Simone Naldi\thanks{Technische Universit\"at Dortmund, Fakult\"at f\"ur Mathematik, Dortmund, Germany} \and
Mohab Safey El Din\thanks{Sorbonne Universit\'es, UPMC Univ Paris 06, CNRS, INRIA Paris Center, LIP6, Equipe PolSys, F-75005, Paris, France.}}

\date{\today}

\maketitle

\begin{abstract}
\noindent
This document describes our freely distributed Maple library {\sc spectra}, for Semidefinite Programming solved Exactly with Computational Tools of Real Algebra. It solves linear matrix inequalities with symbolic computation in exact arithmetic and it is targeted to small-size, possibly degenerate problems for which symbolic infeasibility or feasibility certificates are required. 

\begin{center}
{\bf Keywords}
\end{center}
Computer algebra, symbolic computation, linear matrix inequalities, semidefinite programming, low rank matrices, real algebraic geometry.
\end{abstract}

\newpage

\section{Introduction}

\begin{figure}[ht!]
\centerline{\includegraphics[width=0.7\textwidth]{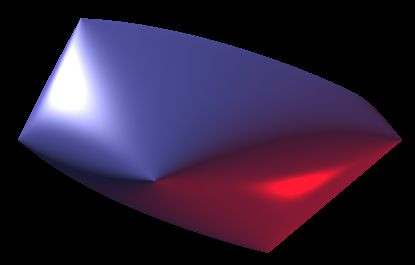}}
\caption{A spectrahedron.\label{fig:spectrahedron}}
\end{figure}

Given symmetric matrices $A_0, A_1, \ldots, A_n$ of size $m$ with rational coefficients, let
\[
{\mathscr S} := \{x \in {\mathbb R}^n : A(x) := A_0 + A_1 x_1 + \cdots + A_n x_n \succeq 0\}
\]
denote the corresponding convex {\em spectrahedron}, defined by the linear matrix inequality (LMI) enforcing that $A$ is positive semidefinite, or equivalently that the eigenvalues of $A$, as functions of $x$, are all nonnegative. Spectrahedra are a broad generalization of polyhedra \cite{v14}. Like polyhedra, spectrahedra have facets, edges and vertices. However, while the facets of a polyhedron are necessarily flat, the facets of a spectrahedron can be curved outwards or inflated, see Figure \ref{fig:spectrahedron} for an example of a spectrahedron of dimension $n=3$ defined by an LMI of size $m=5$. 

Optimization of a linear function on a spectrahedron is called semidefinite programming (SDP), a broad generalization of linear programming (LP) with many applications in control engineering, signal processing, combinatorial optimization, mechanical structure design, etc, see \cite{vb96,t01}. The algebra and geometry of spectrahedra is an active area of study in real algebraic geometry, especially in connection with the problem of moments and the decomposition of real multivariate polynomials as sums-of-squares (SOS), see \cite{l10,bpt12} and references therein.

Our software {\sc spectra} aims at either proving that ${\mathscr S}$ is empty, or finding at least one point in ${\mathscr S}$, using {\em exact arithmetic}. Contrary to numerical algorithms which are based on approximate computations and floating point arithmetic -- such as the projection and rounding heuristics of e.g. \cite{pp08} or \cite{klyz12} -- {\sc spectra} is exclusively based on computations with exact arithmetic. Since exact computations are potentially expensive, {\sc spectra}  should be used when the number $n$ of variables or the size $m$ of the LMI are small. It should not be considered as a competitor to numerical algorithms such as interior-point methods for SDP. It should be primarily used either in potentially degenerate situations, for example when it is expected that $\mathscr S$ has empty interior, or when a rigorous certificate of infeasibility or feasibility is required. 

The input provided to {\sc spectra} is the set of matrices $A_0, A_1, \ldots, A_n$ with rational coefficients describing the pencil $A$ and hence the spectrahedron $\mathscr S$. If $\mathscr S$ is empty,  {\sc spectra} returns the empty list. Otherwise,
the output generated by {\sc spectra} is a finite set 
\begin{equation}\label{par}
{\mathscr Z} := \left\{\left(\frac{q_1(z)}{q_0(z)}, \frac{q_2(z)}{q_0(z)}, \cdots, \frac{q_n(z)}{q_0(z)}\right) : q(z) = 0, \: z\in{\mathbb C}\right\}.
\end{equation}
described by a collection of univariate polynomials with integer coefficients $q, q_0, q_1,\ldots, q_n \in {\mathbb Z}[z]$ and such that 
$\mathscr Z$ meets $\mathscr S$ in at least one real point $x^*$. Such a description is called a rational parametrization. It allows  to isolate the (generally irrational) coordinates of $x^*$ in rational intervals of length given a priori, as small as required.

If $\mathscr S$ is not empty, {\sc spectra}  is guaranteed to compute a point $x^*$ minimizing the rank of $A$ in $\mathscr S$. It solves exactly the (non-convex) optimization problem
\begin{equation}\label{minrank}
\begin{array}{rcll}
r(A) & := & \min & \mathrm{rank}\:A(x) \\
&& \mathrm{s.t.} & x \in \mathscr S.
\end{array}
\end{equation}
This is in sharp contrast with interior-point methods which are designed to compute a point of maximal rank.

The outline of the remainder of the paper is as follows. In Section \ref{background} we survey some background material and extract the essential theoretical results of \cite{hns15} on which {\sc spectra} relies. In Section \ref{started} we provide instructions to download and install {\sc spectra}, and we illustrate its use on two elementary examples. More advanced examples are described in Section \ref{examples}. The performance of {\sc spectra} on larger examples is reported in Section \ref{performance}. Finally, in Section \ref{syntax} we describe formally the exact input and output syntaxes of {\tt SolveLMI}, the main function of {\sc spectra}. 

\section{Background material and main theoretical results}\label{background}

The algorithm implemented in {\sc spectra} computes points in the {\em determinantal varieties}
\[
{\mathscr D}_r := \{x \in {\mathbb C}^n : \mathrm{rank}\:A(x) \leq r\}
\]
for $r=0,1,\ldots,m-1$. By construction it holds
\[
{\mathscr D}_0 \subset {\mathscr D}_1  \subset \cdots {\mathscr D}_{m-1}.
\]
Since the determinant of $A$ vanishes on the boundary ${\partial \mathscr S}$ of $\mathscr S$, it  holds
\[
{\partial \mathscr S} \subset {\mathscr D}_{m-1} \cap {\mathbb R}^n.
\]
When ${\mathscr S}$ is not empty, the value $r(A)$ of the optimization problem (\ref{minrank}) is the minimum integer $r$ such that ${\mathscr D}_r \cap {\mathbb R}^n$ intersects $\mathscr S$, namely the smallest rank on the spectrahedron. Our main geometrical result \cite[Theorem 2]{hns15} states that the spectrahedron contains at least one of the connected components of the real part of the determinantal variety of smallest rank:
\begin{theorem}[{\bf Smallest rank on a spectrahedron}]\label{th:minrank}
Assume that ${\mathscr S}$ is not empty. Let $\mathscr C$ be a connected component of ${\mathscr D}_{r(A)} \cap {\mathbb R}^n$ such that the intersection ${\mathscr C} \cap {\mathscr S}$ is not empty. Then ${\mathscr C} \subset {\mathscr S}$ and hence ${\mathscr C} \subset ({\mathscr D}_{r(A)} \backslash {\mathscr D}_{r(A)-1}) \cap {\mathbb R}^n$.
\end{theorem}
As a consequence of this result, an algorithm computing at least one point in each connected component of ${\mathscr D}_{r(A)} \cap {\mathbb R}^n$ will compute at least one point in the spectrahedron $\mathscr S$. Since the value of $r(A)$ is not known beforehand in general,  {\sc spectra} proceeds iteratively by computing at least one point in the real determinantal variety ${\mathscr D}_r \cap {\mathbb R}^n$ for increasing values of the expected rank $r=0,1,\ldots,m-1$.

More specifically, {\sc spectra} computes points in the determinantal varieties ${\mathscr D}_r$
by projecting onto the subspace of $x$ variables the {\em incidence varieties}
\[
{\mathscr V}_r:=\{(x,y) \in  {\mathbb C}^n\times {\mathbb C}^{m(m-r)} : A(x)Y(y) = 0, \:\:\mathrm{rank}\:Y(y) = m-r\}
\]
for $r=0,1,\ldots,m-1$. The reader familiar with SDP duality will recognize the classical complementarity conditions, see e.g. \cite{vb96,t01}. The dual matrix
\[
Y(y)=
\left(
\begin{array}{ccc}
y_{1,1} & \hspace{-0.2cm} \cdots \hspace{-0.2cm} & y_{1,m-r} \\
\vdots & \hspace{-0.2cm} \hspace{-0.2cm}        & \vdots  \\
y_{m,1} & \hspace{-0.2cm} \cdots \hspace{-0.2cm} & y_{m,m-r}
\end{array}
\right)
\]
depends linearly on the dual variables $y$, and some normalization constraint should be added to ensure that $\mathrm{rank}\:Y(y) = m-r$. Unlike ${\mathscr D}_r$, the incidence variety ${\mathscr V}_r$, up to genericity conditions on the pencil $A$ , turns
to be {\sl smooth and equidimensional}. This crucial geometric property allows for the application of a recursive method which is guaranteed to find at least one point in each connected component of the incidence variety. This leads to the main algorithmic result \cite[Theorem 3]{hns15} on which {\sc spectra} relies:
\begin{theorem}[\bf Exact algorithm for finding a point in a spectrahedron]
Suppose that for each $r=0,1,\ldots,m-1$, the incidence variety ${\mathscr V}_r$ is
smooth and equidimensional and that its defining polynomial system
generates a radical ideal. Suppose also that the determinantal variety ${\mathscr D}_r$ has the expected dimension $n-\binom{m-r+1}{2}$. Then, there is a probabilistic algorithm that takes $A$ as input and returns:
\begin{enumerate}
\item either the empty list, if and only if $\mathscr S$ is empty, or
\item
  a vector $x^*$ such that $A(x^*)=0$, if and only if the linear system $A(x)=0$ has
  a solution, or
\item
   a rational parametrization $q,q_0,q_1,\ldots,q_n \in {\mathbb Z}[z]$
   such that there exists $z^* \in {\mathbb R}$ with $q(z^*)=0$ and: 
\begin{itemize}
  \item $A(q_1(z^*)/q_0(z^*),\ldots,q_n(z^*)/q_0(z^*)) \succeq 0$ and
  \item $\mathrm{rank}\:A(q_1(z^*)/q_0(z^*),\ldots,q_n(z^*)/q_0(z^*)) = r(A)$.
  \end{itemize}
\end{enumerate}
\end{theorem}
The probabilistic nature of the algorithm comes from random changes of variables performed during
the procedure, allowing to put the sets ${\mathscr D}_r$ in generic position.

Recall that the incidence varieties ${\mathscr V}_r$ are defined by enforcing a full column rank constraint on the dual matrix $Y(y)$.
In {\sc spectra} this is achieved as follows \cite[Section 3.1]{hns15}: given a subset of $m-r$ dinstinct integers between $1$ and $r$, we enforce the submatrix of $Y(y)$ whose rows are indexed by these integers to be equal to the identity matrix of size $m-r$. For a given value of $r$, there are $\binom{m}{r}$ distinct choices of row indices and hence the same number of normalized incidence varieties. For each value of $r$, the algorithm in {\sc spectra} processes iteratively these normalized incidence varieties.

Finally, let us explain briefly how {\sc spectra} is able to certify the correctness of the output. This explanation was not included in our paper \cite{hns15}, but we believe it is useful for readers interested in the implementation details. For each computed solution $(x^*, y^*)$ belonging to a connected component of an incidence variety, {\sc spectra} uses exact arithmetic to decide whether $A(x^*)$ is positive semidefinite and to evaluate the rank of $A(x^*)$. If $A(x^*)$ is not positive semidefinite, then the point $x^*$ is discarded. From Theorem \ref{th:minrank} we know that at least one computed point $x^*$ lies on the spectrahedron $\mathscr S$, and this point is of minimal rank, i.e. it solves problem (\ref{minrank}). 

We first build the following characteristic polynomial:
\[
s \mapsto p(s,x) = \det(s I_m+A(x)) = s^m + p_1(x)s^{m-1} + \cdots + p_{m-1}(x)s + p_m(x),
\]
where $I_m$ is the identity matrix of size $m$.
The coefficient $p_k(x) \in \mathbb{Q}[x]$ has degree $k$ in $x$ and it is the $k$-th elementary
symmetric polynomial of the eigenvalues of $A(x)$, for $k=1, \ldots, m$. For example, $p_1(x)$ is the trace of $A(x)$ and $p_m(x)$ is the determinant of $A(x)$. This computation is done only once.

Let $x^* \in \mathbb{R}^n$ be given. The rank defect of $A(x^*)$ is equal to the number of consecutive zeros in the sequence $p_m(x^*), p_{m-1}(x^*), \ldots$ Moreover, by Descartes' rule of signs, $A(x^*) \succeq 0$ if and only if $p_k(x^*) \geq 0$ for all $k=1,\ldots,m$. Hence, computing exactly the rank of $A(x^*)$ and checking its positive semidefiniteness amounts to
determining the signs of $p_k(x^*)$ for $k=1, \ldots, m$.

Whereas this sign determination is a delicate issue when using floating arithmetic and approximate computation, it can be done exactly with {\sc spectra}, since  we represent the point $x^*$ with a rational univariate parametrization with coefficients in $\mathbb{Z}$. 
Indeed, suppose that $x^*$ belongs to the finite set ${\mathscr Z}$ described as in (\ref{par}) by the integer coefficient polynomials $q,q_0,q_1,\ldots,q_n$. Together with the rational intervals isolating each entry of $x^*$, {\sc spectra} computes rational intervals isolating each coefficient $p_k(x^*)$.  Each isolating interval is gradually reduced, until it is so small that at the interval bounds the coefficient takes 1) distinct signs, in which case it vanishes somewhere in the interval, or 2) the same sign, in which case it does not vanish in the whole interval.

\section{Getting started}\label{started}

{\sc spectra}  is freely available as a library for Maple version 16 and above. It can be downloaded in the form of single binary file {\tt SPECTRA.mla} from the following page

\begin{center}
{\tt homepages.laas.fr/henrion/software/spectra}
\end{center}

{\sc spectra} relies on {\sc FGb}, a library for fast computation of Gr\"obner bases, whose Maple interface must be installed, see \cite{f10}.
{\sc spectra} does not work without {\sc FGb}.

In a Maple worksheet, from the directory containing the file {\tt SPECTRA.mla}, please type the command
\begin{verbatim}
> with(SPECTRA);
\end{verbatim}
to activate the main function {\tt SolveLMI}.

\subsection{Half disk}

Let
\[
A(x) = \left(\begin{array}{ccc}
1+x_1 & x_2 & 0 \\
x_2 & 1-x_1 & 0 \\
0 & 0 & x_1
\end{array}\right)
\]
with $n=2$ and $m=3$.
The corresponding spectrahedron
\[
{\mathscr S} = \{ x \in {\mathbb R}^2 : A(x) \succeq 0\} = \{ x \in {\mathbb R}^2 : 1-x^2_1-x^2_2 \geq 0, \: x_1 \geq 0\}
\]
is a half disk. To find a point in $\mathscr S$, we use {\sc spectra} as follows:
\begin{verbatim}
> A := Matrix([[1+x1, x2, 0], [x2, 1-x1, 0], [0, 0, x1]]):
> SolveLMI(A);              
  [[x1 = [0, 0], x2 = [1, 1]]]
\end{verbatim}
This returns the point
\[
x=[0,1] \in {\mathscr S}
\]
in interval arithmetic notation, i.e. 
\[
 x_1 \in [0,0], \: x_2 \in [1,1]
\]
and for each component in $x$ we obtain rational (exact) lower and upper bounds. Here the bounds coincide as the point has rational coordinates.

At this point, matrix $A(x)$ is guaranteed to have minimal rank over all points in ${\mathscr S}$. This rank can be obtained as follows:
\begin{verbatim}
> SolveLMI(A,{rnk});              
  [[x1 = [0, 0], x2 = [1, 1], rnk = 1]]
\end{verbatim}

\subsection{Degenerate spectrahedra}

Let us modify the bottom right entry in the matrix of the previous section, so that now
\[
A(x) = \left(\begin{array}{ccc}
1+x_1 & x_2 & 0 \\
x_2 & 1-x_1 & 0 \\
0 & 0 & x_1-1
\end{array}\right)
\]
and the corresponding spectrahedron ${\mathscr S} =  \{x \in {\mathbb R}^2 : A(x) \succeq 0\} = \{[1,0]\}$ reduces to a point in the plane.
{\sc spectra} can easily deal with such a degenerate case:
\begin{verbatim}
> A := Matrix([[1+x1, x2, 0], [x2, 1-x1, 0], [0, 0, x1-1]]):
> SolveLMI(A);
  [[x1 = [1, 1], x2 = [0, 0]]]
\end{verbatim}

Now let us modify further the bottom right entry, letting
\[
A(x) = \left(\begin{array}{ccc}
1+x_1 & x_2 & 0 \\
x_2 & 1-x_1 & 0 \\
0 & 0 & x_1-2
\end{array}\right)
\]
so that the corresponding spectrahedron is empty. {\sc spectra} returns the empty list, and this is a certificate of emptiness:
\begin{verbatim}
> A := Matrix([[1+x1, x2, 0], [x2, 1-x1, 0], [0, 0, x1-2]]):
> SolveLMI(A);
  []
\end{verbatim}

Since {\sc spectra} is based on exact arithmetic, it is not sensitive to numerical roundoff errors or small parameter changes:
\begin{verbatim}
> A := Matrix([[1+x1, x2, 0], [x2, 1-x1, 0], [0, 0, x1-1-10^(-20)]]):
> SolveLMI(A);
  []
> A:=Matrix([[1+x1, x2, 0], [x2, 1-x1, 0], [0, 0, x1-1+10^(-20)]]):
> SolveLMI(A);
  [[x1 = [36893488147418995335 / 36893488147419103232,
          4611686018427401391 / 4611686018427387904],
    x2 = [-350142318592414077 / 2475880078570760549798248448,
          -2801138548739304423 / 19807040628566084398385987584]]
\end{verbatim}
Displayed with 10 significant digits, the latter point reads:
\[
\begin{array}{l}
x_1 \in [\frac{3689348814741899533}{36893488147419103232}, 
\frac{4611686018427401391}{4611686018427387904}] \approx 1.000000000, \\[1em]
x_2 \in [\frac{-350142318592414077}{2475880078570760549798248448},
\frac{-2801138548739304423}{19807040628566084398385987584}]
\approx -0.1414213562\cdot 10^{-9}.
\end{array}
\] 
The above point is an irrational solution, and the rational intervals are provided so that their floating point approximations are correct up to the number of digits specified in the Maple environment variable {\tt Digits}, which is by default equal to 10. Use the command
\begin{verbatim}
> Digits:=100:
\end{verbatim}
prior to calling {\tt SolveLMI} if you want an approximation correct to 100 digits. At the price of increased computational burden, {\sc spectra} then provides larger integer numerators and denominators in the coordinate intervals.

\section{Examples}\label{examples}

\subsection{Irrational spectrahedron}

In general, each coordinate of a point computed by {\sc spectra} is an algebraic number, i.e. the root of a univariate polynomial with integer coefficients. 

For the classical univariate matrix
\[
A(x_1) =
\begin{bmatrix}
1 & x_1 & 0 & 0 \\
x_1 & 2 & 0 & 0 \\
0 & 0 & 2x_1 & 2 \\
0 & 0 & 2 & x_1
\end{bmatrix}
\]
the spectrahedron reduces to the irrational point $x_1=\sqrt{2}$.
The simple call
\begin{verbatim}
> A:=Matrix([[1, x1, 0, 0], [x1, 2, 0, 0], [0, 0, 2*x1, 2], [0, 0, 2, x1]]):
> SolveLMI(A);            
  [[x1 = [26087635650665550353 / 18446744073709551616,
          13043817825332807945 / 9223372036854775808]]]          
\end{verbatim}
returns an interval enclosure valid to 10 digits.
We can however obtain an exact representation of this point via a rational parametrization:
\begin{verbatim}
> SolveLMI(A, {par});
  [[x1 = [..], par = [_Z^2-2,_Z,[2]]]]
\end{verbatim}
The output parameter {\tt par} contains three univariate polynomials $q, q_0, q_1$ such that the computed point is contained in the finite set
\[
{\mathscr Z}=\{q_1(z)/q_0(z) : q(z) = 0\}=\{2/z : z^2-2 = 0\}=\{\pm\sqrt 2\}
\]
as in (\ref{par}).
Here obviously the rational interval isolates the irrational point $x_1=\sqrt 2$.

\subsection{Algebraic degree}

The algebraic degree of semidefinite programming was studied in \cite{nrs10}. Let us consider the spectrahedron of Example 4 in this reference, for which
\[
A(x) =
\begin{bmatrix}
1+x_3 & x_1+x_2 & x_2 & x_2+x_3 \\
x_1+x_2 & 1-x_1 & x_2-x_3 & x_2 \\
x_2 & x_2-x_3 & 1+x_2 & x_1+x_3 \\
x_2+x_3 & x_2 & x_1+x_3 & 1-x_3
\end{bmatrix}
\]
The following point can be easily found with {\sc spectra}, and it has rank 2, which is guaranteed to be the minimal rank achieved in the spectrahedron:
\begin{verbatim}
> A:=Matrix([[1+x3, x1+x2, x2, x2+x3], [x1+x2, 1-x1, x2-x3, x2],
             [x2, x2-x3, 1+x2, x1+x3], [x2+x3, x2, x1+x3, 1-x3]]): 
> SolveLMI(A, {rnk});                                                                                
  [[x1 = [29909558235590963953/36893488147419103232,
          29909558235593946897/36893488147419103232], 
    x2 = [-18555206088021567643/36893488147419103232,
          -9277603044010395249/18446744073709551616],     
    x3 = [-12556837519724045701/36893488147419103232,
          -12556837519723709525/36893488147419103232],
    rnk = 2]]
\end{verbatim}
With the following instruction we can indeed certify that there is no point of rank 1 or less:
\begin{verbatim}
> SolveLMI(A, {}, [1]);                                                                                
  []
\end{verbatim}

The command
\begin{verbatim}
> SolveLMI(A, {par});
\end{verbatim}
returns the following rational univariate parametrization (\ref{par}) of the above rank 2 point: 
\[\small
\begin{array}{l}
q(z)=16144 z^{10}+ 35160 z^9  + 14536 z^8 - 17690  z^7 - 16278 z^6- 2001 z^5 + 1556 z^4 + 454 z^3 + 23 z^2 -4 z - 1 \\
q_0(z)=161440 z^9 + 316440 z^8  + 116288 z^7 - 123830 z^6  - 97668 z^5  - 10005 z^4  + 6224 z^3  + 1362 z^2  + 46 z - 4 \\
q_1(z) = 97248 z^9  + 146144 z^8 - 18192 z^7  - 134826 z^6  - 63302 z^5  + 4048 z^4 + 6758 z^3  + 846 z^2  - 49 z - 14\\
q_2(z) = 34456 z^9  + 37516 z^8  - 8734 z^7 - 22150 z^6  - 8223 z^5  - 3978 z^4  - 1324 z^3  + 104 z^2  + 103 z + 13\\
q_3(z) =-35160 z^9  - 29072 z^8  + 53070 z^7  + 65112 z^6  + 10005 z^5  - 9336 z^4 - 3178 z^3  - 184 z^2  + 36 z + 10
\end{array}
\]
The degree of the polynomial $q$ in this parametrization can be obtained with the command
\begin{verbatim}
> SolveLMI(A, {deg});
\end{verbatim}

We can obtain more points in the spectrahedron as follows: 
\begin{verbatim}
> SolveLMI(A, {all, rnk, deg}, [2]);
\end{verbatim}
This returns 4 feasible solutions of rank $r=2$, all parametrized by the above degree 10 polynomial.
Notice that this degree matches with the algebraic degree of a generic semidefinite programming problem with parameters
$(m,n,r)=(4,3,2)$, which is 10 according to \cite[Table 2]{nrs10}.  

\subsection{Reproducibility}

Consider the matrix
\[
A(x) = \left(\begin{array}{cc}
1+x_1 & x_2 \\ x_2 & 1-x_1
\end{array}\right)
\]
modeling the unit disk. Two consecutive calls to {\tt SolveLMI} return two distinct points:
{\scriptsize
\begin{verbatim}
> A:=Matrix([[1+x1,x2],[x2,1-x1]]):
> SolveLMI(A); 
 [[x1 = [-21201056044062027875/36893488147419103232, -662533001376936933/1152921504606846976],
   x2 = [-7548363607018988253/9223372036854775808, -1887090901754742967/2305843009213693952]]]
> SolveLMI(A);
   [[x1 = [-10862500438565607907/590295810358705651712, -21725000877131177215/1180591620717411303424],
     x2 = [-576363141759828805/576460752303423488, -9221810268157244495/9223372036854775808]]]
\end{verbatim}
}
After another call, or on your own computer, these intervals should still differ as {\sc spectra} makes random changes of coordinates to ensure that the geometric objects computed are in general position. This kind of behavior is expected when there are infinitely many points of minimal rank in the spectrahedron.

To generate reproducible outputs, the instruction {\tt randomize} can be used to seed the random number generator used by Maple:
{\scriptsize
\begin{verbatim}
> randomize(31415926):
> SolveLMI(A); 
  [[x1 = [-35204733513421104993/36893488147419103232, -35204733513421000447/36893488147419103232],
    x2 = [-2758579864857623899/9223372036854775808, -5517159729715231413/18446744073709551616]]]
> randomize(31415926):
> SolveLMI(A); 
  [[x1 = [-35204733513421104993/36893488147419103232, -35204733513421000447/36893488147419103232],
    x2 = [-2758579864857623899/9223372036854775808, -5517159729715231413/18446744073709551616]]]
\end{verbatim}
}

\subsection{Convex quartic}

\begin{figure}[!ht]
\centering
\includegraphics[width=0.5\textwidth]{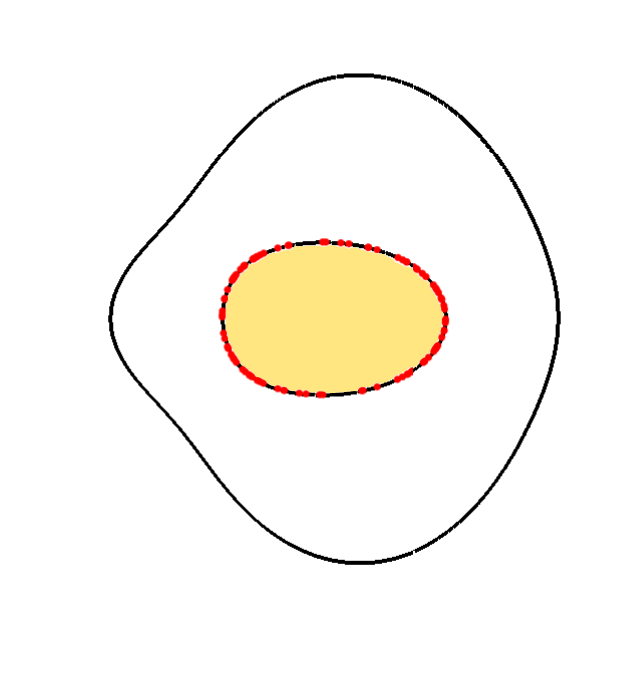}
\caption{Quartic curve (black) with sample points (red) on the boundary of its spectrahedron (orange).\label{fig:quartic}}
\end{figure}

Let
\[
A(x) = \left(\begin{array}{cccc}
1+x_1& x_2 & 0 & 0 \\
x_2 & 1-x_1 & x_2 & 0 \\
0 & x_2 & 2+x_1 & x_2 \\
0 & 0 & x_2 & 2-x_1
\end{array}\right).
\]
The spectrahedron ${\mathscr S} = \{x \in {\mathbb R}^2 : A(x) \succeq 0\}$ is the orange region whose boundary is the internal oval of the smooth quartic determinantal curve $\{x \in {\mathbb R}^2 : \mathrm{det}\:A(x) = 0\}$ represented in black on Figure \ref{fig:quartic}.
With the following instructions
\begin{verbatim}
> A:=Matrix([[1+x1,x2,0,0],[x2,1-x1,x2,0],[0,x2,2+x1,x2],[0,0,x2,2-x1]]):
> SolveLMI(A,{},[3]);
> SolveLMI(A,{},[3]);
> ...
\end{verbatim}
we compute several points on the boundary of $\mathscr S$, they are plotted in red on Figure \ref{fig:quartic}.
Note the third input argument {\tt [3]} which specifies to {\tt SolveLMI} the expected rank of the computed point.
Since the determinantal curve is smooth, we know that the rank of $A(x)$ equals
3 on the whole curve, and in particular on the boundary of $\mathscr S$.
Since the rank is specified, {\sc spectra} does not have to process iteratively the incidence varieties
corresponding to points of smaller ranks, thereby reducing the computational burden to find at least one point in the spectrahedron.

Each of these points is represented by a rational univariate parametrization of degree $12$, obtained with the instruction
\begin{verbatim}
> SolveLMI(A,{par},[3]);
\end{verbatim}
For example,
for the point $(x_1,x_2) \approx (-0.9689884394,\allowbreak -0.2434013983)$ the polynomial $q$ in the rational
parametrization (\ref{par}) is
\begin{align*}
q(z) & = 5506034827600\,z^{12} - 4608031295324\,z^{10} - 192908794368\,z^{9} + 25693318717857\,z^{8} + \\
     & + 4774492660608\,z^{7} - 17188212283956\,z^{6} - 23438418515712\,z^{5} + 64967482316484\,z^{4} - \\
     & - 11285164470528\,z^{3} - 11887630039728\,z^{2} + 296990121024.
\end{align*}
Recall that the algebraic degree of a point $x^*$ in $\mathscr S$
is the degree of the minimal algebraic extension of the ground field (here the rational numbers) required to represent $x^*$. The algebraic degree depends on the size of the pencil $A$ but also on the rank $r$ of $A(x^*)$. With $(m,n,r)=(4,2,3)$ and generic data, the algebraic degree is $12$, cf. \cite[Table 2]{nrs10}, which indeed coincides with the degree of the exact representation of $x^*$ 
computed  by {\sc spectra}.

\subsection{Polynomial sums of squares}

Deciding whether a multivariate real polynomial is non-negative is difficult in general. A sufficient condition, or certificate for non-negativity, is that the polynomial can be expressed as a sum of squares (SOS) of other polynomials. Finding a polynomial SOS decomposition amounts to finding a point in a specific spectrahedron called Gram spectrahedron, see e.g. \cite{cpsv16} and references therein.

As an example, consider the homogeneous ternary quartic
\[
f(u) = u_1^4+u_1u_2^3+u_2^4-3u_1^2u_2u_3-4u_1u_2^2u_3+2u_1^2u_3^2+u_1u_3^3+u_2u_3^3+u_3^4.
\]
The polynomial $f$ belongs to a series of examples provided by C. Scheiderer in \cite{s16} to answer (in the negative) the following question by B. Sturmfels: let $f$ be a polynomial with rational coefficients which is an SOS of polynomials with real coefficients; is it an SOS of polynomials with rational coefficients?
Scheiderer's counterexamples prove that, generally speaking, there is no hope in obtaining nonnegativity
certificates over the rationals. However, certificates exist in some algebraic extension of the field of rational numbers.

In the graded reverse lexicographic ordered monomial basis,
the Gram matrix of $f$ is the matrix
\[
A(x)=
\left(
\begin{array}{cccccc}
1 & 0 & x_1 & 0 & -{3}/{2}-x_2 & x_3 \\
0 & -2x_1 & {1}/{2} & x_{2} & -2-x_4 & -x_5 \\
x_1 & {1}/{2} & 1 & x_4 & 0 & x_6 \\
0 & x_2 & x_4 & -2x_3+2 & x_5 & {1}/{2} \\
-{3}/{2}-x_2 & -2-x_4 & 0 & x_5 & -2x_6 & {1}/{2} \\
x_3 & -x_5 & x_6 & {1}/{2} & {1}/{2} & 1
\end{array}
\right)
\]
depending linearly on 6 real parameters.
The Gram spectrahedron ${\mathscr S} = \{x \in {\mathbb R}^6 : A(x) \succeq 0\}$ parametrizes the set of all SOS decompositions of $f$. We deduce by the discussion
above that $\mathscr S$ does not contain rational points. In particular, its interior is empty.

Let us use {\sc spectra} to compute points in $\mathscr S$ and hence to get positivity certificates for $f$:
{\scriptsize
\begin{verbatim}
> A := Matrix([[1,0,x1,0,-3/2-x2,x3], [0,-2*x1,1/2,x2,-2-x4,-x5], [x1,1/2,1,x4,0,x6],
               [0,x2,x4,-2*x3+2,x5,1/2], [-3/2-x2,-2-x4,0,x5,-2*x6,1/2], [x3,-x5,x6,1/2,1/2,1]]):
> SolveLMI(A, {rnk, deg, par});
[[[x1 = [..], x2 = [..], x3 =  [..], x4 = [..], x5 = [..], x6 = [..]],
    rnk = 2, deg = 3,
    par = [8*z^3+8*z+1, 24*z^2-8, [16*z+3, -24*z^2+8, 8*z^2+6*z+8, -16*z^2+6*z+16, -16*z-3, 16*z+3]]]
\end{verbatim}
}
We obtain an irrational point $x \in {\mathscr S}$ whose coordinates are algebraic numbers of degree 3, and which belongs to the finite set
\[
\left\{\left(\frac{16z+3}{24z^2-8},\frac{-24z^2+8}{24z^2-8},\frac{8z^2+6z+8}{24z^2-8},\frac{-16z^2+6z+16}{24z^2-8},\frac{-16z-3}{24z^2-8},\frac{16z+3}{24z^2-8}\right) : 8z^3-8z-1=0\right\}
\]
At this point, the Gram matrix $A$ has rank 2, and hence $f$ is an SOS of 2 polynomials.

Let us compute more non-negativity certificates of rank 2:
\begin{verbatim}
> SolveLMI(A,{rnk,deg,par,all},[2]);
\end{verbatim}
In addition to the point already obtained above, we get another point. The user can compare this output
with \cite[Ex.\,2.8]{s16}: it turns out that these are the only 2 points of rank 2. Other points in 
the Gram spectrahedron have rank 4 and they are convex combinations of these 2.

\section{Performance}\label{performance}

\subsection{Exponential bit-size spectrahedron}

For a given $n \in {\mathbb N}$, consider the spectrahedron
\[
{\mathscr S}_n = \left\{x \in {\mathbb R}^n :
\left(\begin{array}{cc}
1 & 2 \\
2 & x_1
\end{array}\right) \succeq 0,
\left(\begin{array}{cc}
1 & x_1 \\
x_1 & x_2
\end{array}\right) \succeq 0,
\cdots, 
\left(\begin{array}{cc}
1 & x_{n-1} \\
x_{n-1} & x_n
\end{array}\right)\succeq 0
\right\}.
\]
For every $x^* \in {\mathscr S}_n$ it holds
$
x^*_n \geq (x^*_{n-1})^2 \geq \cdots \geq (x^*_1)^{2^{n-1}} \geq 2^{2^n},
$
which shows that exponentially many bits are required to represent a point.
It is elementary to check that each of the above $n$ matrices of size 2 can have rank 1,
and hence that we can compute a point of rank $n$ as follows:
\begin{verbatim}
> with(LinearAlgebra):
> A:=DiagonalMatrix([<<1,2>|<2,x1>>,<<1,x1>|<x1,x2>>,<<1,x2>|<x2,x3>>,..]):
> SolveLMI(A,{},[n]);
\end{verbatim}
Recall from Section \ref{background} that {\sc spectra} examines iteratively a family of $\binom{m}{r}= \binom{2n}{n}$
incidence varieties, a number growing exponentially with $n$. For example there are $12870=\binom{16}{8}$ incidence
varieties to test to solve our problem for $n=8$. Hence we could expect {\sc spectra} to perform poorly on this example.
However, on our standard desktop PC equipped with Intel i7 processor at 2.5GHz and 16GB RAM, we were able to handle
spectrahedra of size $2n=10$ in 29 seconds, and of size $2n=12$ in 505 seconds.

\subsection{Random spectrahedra}

Finally, we report on randomly generated examples. The rational entries of $A$ are generated as quotients of integers
drawn uniformly in the interval $[-100,100]$. Here is the script we used to generate a random symmetric
pencil given the number $n$ of variables and the size $m$:
\begin{verbatim}
> var:=[seq(cat('x',i),i=1..n)]:
> A:=Matrix(m,m):
> for i from 1 to m do
     for j from i to m do
       A[i,j]:=randpoly(var, degree=1, dense, coeffs=rand(-100..100)):
       A[j,i]:=A[i,j];
     od:
   od:
\end{verbatim}
For each instance, given the expected rank $r$, we execute the command
\begin{verbatim}
> SolveLMI(A,{},[r]);
\end{verbatim}

\begin{figure}[!ht]
\centering
\begin{tikzpicture}
\begin{loglogaxis} [xlabel = number of variables, ylabel = time (seconds)]
\addplot [ultra thick, red]
file {m2r1.txt};
\end{loglogaxis}
\end{tikzpicture}
\caption{Timings for random instances of size $m=2$ and rank $r=1$, as a function of the number of variables $n$.}
\label{fig:m2r1}
\end{figure}
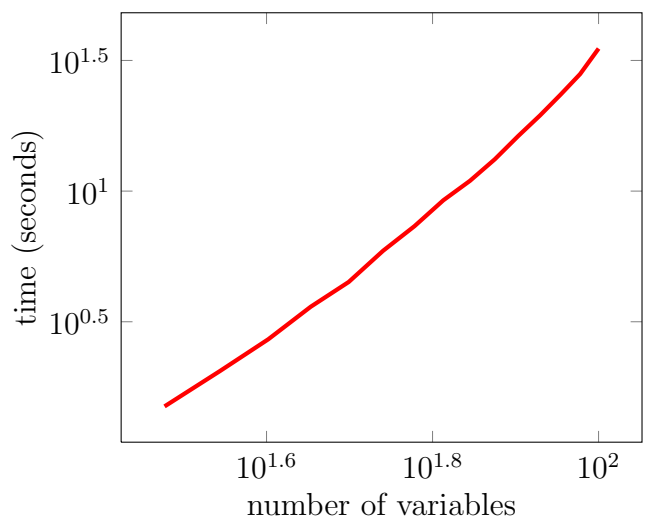

For $m = 2$, $r = 1$ and values of $n$ ranging from 30 to 100, we obtain the timings reported on Figure \ref{fig:m2r1}. 
This corresponds to spectrahedra whose boundaries belong to determinantal varieties of increasing dimension. Moreover,
the singularity locus of the determinant has positive dimension, it is a linear subspace of co-dimension $3$.
We observe a polynomial dependence of the computational time as a function of the number of variables, with exponent around $3$.

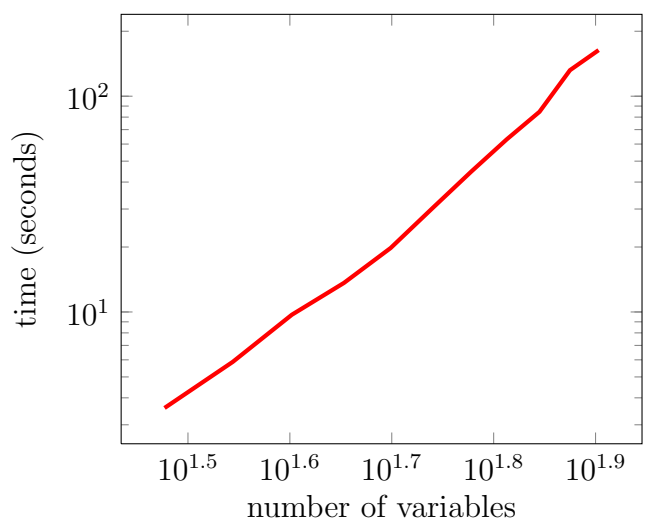
\begin{figure}[!ht]
\centering
\begin{tikzpicture}
\begin{loglogaxis} [xlabel = number of variables, ylabel = time (seconds)]
\addplot [ultra thick, red]
file {m3r2.txt};
\end{loglogaxis}
\end{tikzpicture}
\caption{Timings for random instances of size $m=3$ and rank $r=2$, as a function of the number of variables $n$.}
\label{fig:m3r2}
\end{figure}

When $m = 3$, $r = m-1 = 2$ and values of $n$ ranging from 30 to 80, we obtain the timings reported on Figure \ref{fig:m3r2}, depending polynomially on $n$ with an exponent around $4$.

\section{Input and output syntax}\label{syntax}

\subsection{Input}

The calling sequence of function {\tt SolveLMI} is as follows:
\begin{verbatim}
> SolveLMI(A, options, ranks);
\end{verbatim}
where
\begin{itemize}
\item
$A$ is a symmetric matrix of size $m$ with rational coefficients, depending affinely on $n$ variables; 
\item
{\tt options} (optional) is a set that can contain the following keywords:
\begin{itemize}
\item[{\tt all}]: compute as many solutions as possible, which can be computationally demanding; when this option is not specified, the algorithm is stopped as soon as one solution is computed, which is typically much faster;
\item[{\tt rnk}]: return the rank of $A$ at every computed solution;
\item[{\tt par}]: return the rational univariate parametrization of every computed solution;
\item[{\tt deg}]: return the algebraic degree of every computed solution;
\end{itemize}
\item
{\tt ranks} (optional) is a list of nonnegative integers corresponding to expected ranks of computed solutions. The default value is $[0,1,\ldots,m-1]$. The algorithm is run for each value $r$ in {\tt ranks} by solving the quadratic system of equations
\[A(x)Y(y)=0\]  for a vector $x$ and a matrix $Y(y)$ with $m$ rows and $m-r$ columns whose entries are stored in a vector $y$. It may happen that the rank of $A(x)$ at a computed solution $x$ is strictly less than $r$.

\end{itemize}

\subsection{Output}

Let us denote by $x_1, x_2, \ldots, x_n$ the variables on which matrix $A$ depends affinely. They are gathered in a vector $x \in {\mathbb R}^n$. When the input argument {\tt options} is empty, the output returned by {\tt SolveLMI} is
\begin{itemize}
\item eithter the empty list {\tt []} in which case  $\mathscr S$ is empty, or 
\item a rational enclosure of a single point $x \in {\mathscr S}$, in the form
\begin{verbatim}
> SolveLMI(A)
  [[x1 = [a1, b1], x2 = [a2, b2], ..., xn = [an, bn]]]
\end{verbatim}
where $a_i, b_i$ are rational numbers,  displayed as ratios of integers. This means that each coordinate $x_i$ belongs to the interval $[a_i,b_i]$ ensuring a floating point approximation of $x$ valid to a number of digits equal to the Maple environment variable {\tt Digits}. When $a_i=b_i$ this implies that $x_i$ is a rational number.
\end{itemize}

When {\tt options} contains the keyword {\tt all}, more points can be returned, in the form of a list
\begin{verbatim}
> out := SolveLMI(A, {all})
  [[x1 = [a11, b11], x2 = [a12, b12], ..., xn = [a1n, b1n]],
   [x1 = [a21, b21], x2 = [a22, b22], ..., xn = [a2n, b2n]], 
   ...]
\end{verbatim}
such that {\tt nops(out)} is the number of computed points, {\tt out[1]} is the first point, {\tt out[2]} is the second point, etc.

When {\tt options} contains the keyword {\tt rnk}, the rank of $A$ at $x$ is returned:
\begin{verbatim}
> SolveLMI(A, {rnk})
  [[x1 = [a1, b1], x2 = [a2, b2], ..., xn = [an, bn], rnk = r]]
\end{verbatim}

These keywords and the following ones can be freely combined:
\begin{verbatim}
> SolveLMI(A, {all, rnk})
  [[x1 = [a11, b11], x2 = [a12, b12], ..., xn = [a1n, b1n], rnk = r1],
   [x1 = [a21, b21], x2 = [a22, b22], ..., xn = [a2n, b2n], rnk = r2], 
   ...]
\end{verbatim}

When {\tt options} contains the keyword {\tt par}, a rational univariate parametrization of $x$ is returned:
\begin{verbatim}
> SolveLMI(A, {par})
  [[x1=[a1,b1], x2=[a2,b2], ..., xn=[an,bn], par=[q,q0,[q1,q2,...,qn]]]
\end{verbatim}
This parametrization is such that $q, q_0, q_1, q_2, \ldots, q_n$ are univariate polynomials with integer coefficients such that $x$ belongs to the finite set
\[
\left\{\left(\frac{q_1(z)}{q_0(z)}, \frac{q_2(z)}{q_0(z)}, \cdots, \frac{q_n(z)}{q_0(z)}\right) : q(z) = 0, \: z\in{\mathbb C}\right\}.
\]
The intervals $[a_i,b_i]$ are provide to isolate the computed point from this set of points.

When {\tt options} contains the keyword {\tt deg}, the degree of the polynomial $q$ in the rational univariate parametrization of each computed point $x$ is also returned:
\begin{verbatim}
> SolveLMI(A, {deg})
  [[x1 = [a1, b1], x2 = [a2, b2], ..., xn = [an, bn], deg = d]]
\end{verbatim}

\end{document}